\documentclass[preprint,12pt]{elsarticle}

\usepackage{amssymb}
\usepackage{amsmath}
\usepackage{subcaption}					
\usepackage{mathtools} 					
\usepackage{algorithm} 					
\usepackage{algpseudocode} 				
\usepackage{booktabs}					
\usepackage{enumitem}       			
\usepackage{numprint}					
\usepackage{xcolor}                     
\usepackage{bm} 						
\usepackage{spalign} 					

\npthousandsep{\,}						

\graphicspath{{./figures/}}

\journal{International Journal of Engineering Science}

\newcommand{\bs}[1]{\boldsymbol{#1}}

\newcommand{\dx}[0]{}
\newcommand{\ds}[0]{}

\newcommand{\R}[0]{\mathbb{R}}

\newcommand{\vertiii}[1]{{\left\vert\kern-0.25ex\left\vert\kern-0.25ex\left\vert #1 \right\vert\kern-0.25ex\right\vert\kern-0.25ex\right\vert}}

\DeclarePairedDelimiterX\Set[1]{\lbrace}{\rbrace}%
 {  #1 }

\newtheorem{remark}{Remark}

\begin{document}

\begin{frontmatter}

	\title{Nitsche methods for constrained problems in mechanics\tnotemark[1]}
	 \tnotetext[1]{This work was supported by the Research Council of Finland (Flagship of Advanced Mathematics for Sensing Imaging and Modelling grant 359181) and the Portuguese government through FCT (Funda\c c\~ao para a Ci\^encia e a Tecnologia), I.P., under the project UIDB/04459/2025.}

 \author[mech]{Tom Gustafsson}
 \ead{tom.gustafsson@aalto.fi}

 \author[msa]{Antti Hannukainen}
 \ead{antti.hannukainen@aalto.fi}

 \author[msa]{Vili Kohonen\corref{cor1}}
 \ead{vili.kohonen@aalto.fi}

 \author[ist]{Juha Videman}
 \ead{jvideman@math.tecnico.ulisboa.pt}

 \cortext[cor1]{Corresponding author}

 \affiliation[mech]{organization={Department of Mechanical Engineering, Aalto University},
             addressline={P.O. Box 11100},
             postcode={00076 Aalto},
             city={Espoo},
             country={Finland}}

 \affiliation[msa]{organization={Department of Mathematics and Systems Analysis, Aalto University},
             city={Espoo},
             country={Finland}}

 \affiliation[ist]{organization={CAMGSD/Departamento de Matem\'atica, Instituto Superior T\'ecnico, Universidade de Lisboa},
             addressline={Av. Rovisco Pais 1},
             postcode={1049-001},
             city={Lisbon},
             country={Portugal}}

\begin{abstract}
We present guidelines for deriving new Nitsche Finite Element Methods to enforce equality and inequality constraints that act on the value of the unknown mechanical quantity.
We first formulate the problem
as a stabilized finite element method for
the saddle point formulation where
a Lagrange multiplier enforces
the underlying constraint. 
The Nitsche method is then presented in a general minimization form,
suitable for adding constraints to nonlinear finite element methods
and allowing straightforward computational implementation
with automatic differentation.
This extends the method beyond classical boundary condition enforcement.
To validate these ideas, we present Nitsche formulations for a range of problems in solid mechanics and give numerical evidence of the convergence rates of the Nitsche method. 
\end{abstract}

\begin{highlights}
\item We derive a novel interpretation of the Nitsche Finite Element Method to solve constrained problems in solid mechanics
\item We use the new formulation to solve several contact problems and validate the convergence rates numerically
\item The computational implementation utilizes a general minimization form and automatic differentation handles the variational derivatives
\end{highlights}


\end{frontmatter}

\section{Introduction}

In the 1970s, Nitsche~\cite{nitsche} introduced a method for enforcing the Dirichlet boundary conditions in the weak formulation of the finite element method.
The benefits of the Nitsche method over its alternatives, the penalty method~\cite{courant1943variational} and the mixed Lagrange multiplier method~\cite{babuvska1973finite}, include consistency,  symmetry and positive definiteness, stability, and conditioning; cf.~\cite{gustafsson2025}.
Hence, the ideas of Nitsche have later found applications also in domain decomposition~\cite{becker2003finite}, discontinuous Galerkin methods~\cite{hansbo2002discontinuous},
unfitted finite element and fictitious domain methods~\cite{burman2015cutfem},
Robin boundary conditions~\cite{juntunen2009nitsche}, Kirchhoff plates~\cite{gustafsson2021nitsche},
elastic contact problems~\cite{chouly2013nitsche}, friction problems~\cite{chouly2014adaptation},
obstacle problems~\cite{gustafsson2017mixed}, elastoplastic torsion problems~\cite{chouly2023nitsche}, and
plate contact problems~\cite{gustafsson2017nitsche,fabre2021nitsche}, to name a few.

In modern literature,
the Nitsche method
is often presented
as a consistency correction to the penalty method, i.e.,
to make sure that the exact solution
satisfies the discrete  weak form, by
including additional terms; cf., e.g., ~\cite{benzaken2024constructing, gustafsson2025}.
This
fails to answer to
the obvious follow-up question:
\emph{how to generalize the
	method to other
	problems and constraints?}
Consequently, we take
a different viewpoint
based on the saddle point
formulation and present guidelines to derive
optimal Nitsche methods for both equality and inequality
constraints.
These guiding principles apply
to arbitrary constrained
minimization problems, where an affine constraint acts
on the values of the unknown.
This extends the method beyond 
conventional boundary condition enforcement.
The resulting energy minimization form 
is particularly well-suited for adding
constraints to
nonlinear problems.
We follow
the guidelines in practice
by deriving and implementing  novel
methods for problems such as
inequality
boundary condition for
plates, membrane-solid contact problem,
and plate-plate contact problem.

Our formulation is based on
a reinterpretation of the Nitsche method
by Stenberg~\cite{stenberg1995}.
More precisely,
on the idea that
the Nitsche method corresponds to
a conforming stabilized Lagrange multiplier
method
(in the sense of residual Barbosa--Hughes stabilization~\cite{barbosa1992circumventing}).
For such stabilized methods,
an element-wise elimination of
the Lagrange multiplier will
lead to the equivalent Nitsche method;
see also \cite{chouly2013nitsche, burman2017galerkin, chouly2023finite, burman2023augmented}.
This connection was exploited in the
error analysis in Gustafsson~\cite{gustafsson2018finite} and
consequently it became evident that the
stabilized formulation leads to
a family of methods with
constraints of rather general form.

Let us summarize the prior theoretical findings for  existing Nitsche methods:
(1) stability becomes automatic if correct
residuals of the Lagrange multiplier are included in the formulation; (2)
optimal convergence rates are obtained
if the scaling of the residual
with respect to the
mesh and the material parameters
correctly mimics the continuous
Sobolev norm of the Lagrange multiplier.

In this paper, built on these observations, we  explain how 
to define and use general Nitsche methods,
present examples of their adaptation
to novel physical problems,
and
show numerical evidence of their
convergence.
A more careful numerical analysis, i.e.,
proofs of stability and error estimates
for the general formulation, is
beyond the scope of this work.
The remainder of the paper is organized as follows. In Section 2, we review
the original Nitsche method for the Dirichlet boundary condition $u = g$.  In Section 3, we
demonstrate how the method is
modified to satisfy
the inequality constraint
$u \geq g$.
In Section 4, we describe a
general form of the Nitsche method
and how it can be used to derive 
new methods.
In Section 5, we present some
methods that exist in the literature
and are of this general form.
In Section 6, we derive novel
Nitsche methods from the general form and show numerical
evidence of their performance.

\section{Background: Dirichlet boundary conditions}

For better context, we briefly recall the original method of Nitsche~\cite{nitsche}.
Consider the model problem
\begin{equation}
	\label{eq:orig}
	\begin{aligned}
		-\nabla \cdot \kappa \nabla u & = f, \quad \text{in $\Omega$},          \\
		u                             & = g, \quad \text{on $\partial \Omega$},
	\end{aligned}
\end{equation}
where $f\in L^2(\Omega)$ is the problem-specific load,
$g\in H^{1/2}(\partial\Omega)$ is the Dirichlet boundary condition
and $\kappa > 0$ is a given parameter.
Note that Nitsche~\cite{nitsche} considered
only the case $\kappa=1$, i.e., normalized units.

Problem \eqref{eq:orig} is solved in a finite element space
$V_h \subset H^1(\Omega)$
where $h$ is the mesh parameter.
The classical Nitsche method for the Dirichlet boundary condition $u=g$ reads: find $u_h \in V_h$
such that
\begin{equation}
	\label{eq:nitscheweak}
	\begin{aligned}
		 & \int_\Omega \kappa \nabla u_h \cdot \nabla v_h\dx + \int_{\partial \Omega} \frac{1}{\gamma} u_h v_h \ds- \int_{\partial \Omega} \kappa \frac{\partial u_h}{\partial n} v_h\ds-\int_{\partial \Omega} u_h\,\kappa \frac{\partial v_h}{\partial n}\ds \\
		 & \quad = \int_\Omega fv_h\dx + \int_{\partial \Omega} \frac{1}{\gamma} g v_h \ds -\int_{\partial \Omega} g  \,\kappa\frac{\partial v_h}{\partial n}\ds \qquad \forall v_h \in V_h.
	\end{aligned}
\end{equation}
Here $\gamma$ is an $L^2(\Omega)$ function whose restriction to each element $K$ is the product of a dimensionless stabilization parameter $\alpha > 0$
and the local mesh parameter $h_K$, divided by $\kappa$, i.e.
$$\gamma|_K = \frac{\alpha h_K}{\kappa}.$$

According to Stenberg~\cite{stenberg1995}, the Nitsche method corresponds to using stabilized finite element method to discretize
the saddle point problem: find $(u, \lambda) \in H^1(\Omega) \times H^{-1/2}(\partial \Omega)$ such that
\begin{equation}
	\label{eq:saddlepoint}
	\mathcal{L}(u, \lambda) = \inf_{v \in H^1(\Omega)} \sup_{\mu \in H^{-1/2}(\partial \Omega)} \mathcal{L}(v, \mu).
\end{equation}
The Lagrangian in \eqref{eq:saddlepoint} is defined as
\begin{equation}
	\label{eq:lagrangian}
	\mathcal{L}(u, \lambda) = J(u) - \langle u - g, \lambda \rangle,
\end{equation}
where
\begin{equation}
	J(u) = \frac12 \int_\Omega \kappa \nabla u \cdot \nabla u \dx  - \int_\Omega fu \dx
\end{equation}
is the energy functional of the original problem \eqref{eq:orig} \cite{brezis2011functional},
$\langle \cdot, \cdot \rangle$ denotes the duality pairing between $H^{1/2}(\partial \Omega)$ and $H^{-1/2}(\partial \Omega)$, and $\lambda$
is a Lagrange multiplier.  
We are able to derive the strong formulation of problem \eqref{eq:saddlepoint} to reveal that
$\lambda = \kappa \frac{\partial u}{\partial n}$.

The Nitsche method \eqref{eq:nitscheweak} corresponds to
the stabilized discrete saddle point problem: find $(u_h, \lambda_h) \in V_h \times Q_h$ such that
\begin{equation}
	\label{eq:discstab}
	\mathcal{L}_h(u_h, \lambda_h)  =\inf_{v_h \in V_h} \sup_{\mu_h \in Q_h} \mathcal{L}_h(v_h, \mu_h)
    \end{equation}
    where
    the stabilized Lagrangian $\mathcal{L}_h$ is defined as
    \begin{equation}
		\label{eq:stabilizedlagrangian}
    \mathcal{L}_h(v_h, \mu_h) = \Big\{\mathcal{L}(v_h, \mu_h) -  \int_{\partial \Omega} \frac{\gamma }{2} \Big( \mu_h - \kappa \frac{\partial v_h}{\partial n}\Big)^2 \ds\Big\}.
\end{equation}
The space $Q_h \subset H^{-1/2}(\partial \Omega)$ denotes a set of
elementwise polynomials that are globally discontinuous.
The variational formulation of problem \eqref{eq:discstab} reads: find $(u_h, \lambda_h) \in V_h \times Q_h$ such that
\begin{equation}
	\begin{aligned}
		 & \int_\Omega \kappa \nabla u_h \cdot \nabla v_h \dx - \int_{\partial \Omega} \lambda_h v_h \ds - \int_{\partial \Omega} u_h \mu_h \ds                                                                                         \\
		 & \quad - \int_{\partial \Omega} \gamma  \Big(\lambda_h - \kappa \frac{\partial u_h}{\partial n}\Big)\Big( \mu_h - \kappa \frac{\partial v_h}{\partial n} \Big)\ds = \int_\Omega fv_h \dx - \int_{\partial \Omega} g \mu_h \ds
	\end{aligned}
\end{equation}
for every $(v_h, \mu_h) \in V_h \times Q_h$.
This method is stable for any pair of finite element spaces $V_h$ and $Q_h$ if $0 < \alpha < C_I$ where $C_I$ is a constant of an inverse estimate, cf.~\cite{stenberg1995}.

Following the steps laid out in Stenberg~\cite{stenberg1995}, we can show that an elimination of the Lagrange multiplier element-by-element leads to the formulation \eqref{eq:nitscheweak}.  Assuming that $\kappa$ is an element-wise constant and
choosing $v_h = 0$ leads to
\begin{equation}
	- \int_{\partial \Omega} u_h \mu_h \ds - \int_{\partial \Omega} \gamma \Big(\lambda_h - \kappa \frac{\partial u_h}{\partial n}\Big) \mu_h\ds = - \int_{\partial \Omega} g \mu_h \ds \quad \forall \mu_h \in Q_h.
\end{equation}
Denoting by $\pi_h$ the $L^2$ projection onto $Q_h$, this reads
\begin{equation}
	- \pi_h u_h - \gamma \Big(\lambda_h - \kappa \pi_h \frac{\partial u_h}{\partial n}\Big) = - \pi_h g
\end{equation}
which gives
\begin{equation}
	\lambda_h = \kappa \pi_h \frac{\partial u_h}{\partial n}-\frac{1}{\gamma}(\pi_h u_h - \pi_h g).
\end{equation}
Let $V_h$ correspond to the space of continuous piecewise polynomial functions of degree $p\geq 1$. If $Q_h$ is chosen as discontinuous piecewise polynomials of degree $p$ or greater, then $\pi_h$ reduces to the identity mapping, i.e.
\begin{equation}
	\label{eq:disclagmult}
	\lambda_h = \kappa \frac{\partial u_h}{\partial n}-\frac{1}{\gamma}(u_h - g).
\end{equation}

Now we can eliminate the Lagrange multiplier
from the discrete saddle point problem.
After substituting $\lambda_h$ to the stabilized Lagrangian $\mathcal{L}_h$ we obtain
\begin{align*}
	J(u_h) & - \int_{\partial \Omega} \Big( \kappa \frac{\partial u_h}{\partial n}-\frac{1}{\gamma}(u_h - g)\Big)(u_h - g) \ds                                                           \\
	       & \quad -  \int_{\partial \Omega} \frac{\gamma}{2} \Big( \kappa \frac{\partial u_h}{\partial n}-\frac{1}{\gamma}(u_h - g) - \kappa \frac{\partial u_h}{\partial n}\Big)^2 \ds
\end{align*}
which simplifies to
\begin{align*}
	J(u_h) & - \int_{\partial \Omega} \kappa \frac{\partial u_h}{\partial n}(u_h - g) + \int_{\partial \Omega}\frac{1}{\gamma}(u_h - g)^2 \ds \\
	       & \quad -  \int_{\partial \Omega} \frac{1}{2 \gamma}(u_h - g)^2 \ds
\end{align*}
and finally to
\begin{equation}
	\label{eq:nitschemin}
	\begin{aligned}
		J_h(u_h) = J(u_h) - \int_{\partial \Omega} \kappa \frac{\partial u_h}{\partial n}(u_h - g)\ds + \frac12 \int_{\partial \Omega}\frac{1}{\gamma}(u_h - g)^2 \ds.
	\end{aligned}
\end{equation}
The classical Nitsche method \eqref{eq:nitscheweak} corresponds to
the following minimization problem \cite{nitsche}: find $u_h \in V_h$ such that
\begin{equation}
	J_h(u_h) = \inf_{v_h \in V_h} J_h(v_h).
\end{equation}

Note that by
adding and subtracting $ \int_{\partial \Omega} \frac{\gamma}{2} (\kappa \frac{\partial u_h}{\partial n})^2 \ds$ to complete the square we can rewrite $J_h$ in \eqref{eq:nitschemin} as
\begin{align*}
	J(u_h) + \int_{\partial \Omega} \frac{\gamma}{2} \Big(\kappa  \frac{\partial u_h}{\partial n} - \frac{1}{\gamma}(u_h - g)\Big)^2\ds -  \int_{\partial \Omega} \frac{\gamma}{2} \Big(\kappa \frac{\partial u_h}{\partial n}\Big)^2 \ds.
\end{align*}
We can formally compare the above to the penalty method (e.g., \cite{gustafsson2025})
which is obtained by minimizing the penalized functional:
\begin{equation*}
	J(u_h) + \frac{1}{2\varepsilon} \int_{\partial \Omega} (u_h - g)^2 \ds, \quad \varepsilon > 0.
\end{equation*}
The penalty method is known to approximate the Robin boundary condition
\begin{equation}
	\kappa \frac{\partial u}{\partial n} = \frac{1}{\varepsilon}(u - g),
\end{equation}
which satisfies $u - g \rightarrow 0$ when $\varepsilon \rightarrow 0$.

\section{Inequality constraints: Signorini's problem}

For inequality constraints, we first consider problem \eqref{eq:orig} with the inequality constraint $u \geq g$ on the boundary $\partial\Omega$.
The resulting problem 
\begin{equation}
	\begin{aligned}
		\label{eq:signorini}
		-\nabla \cdot \kappa \nabla u																	  & = f, \quad \text{in $\Omega$},		  \\
		u \geq g,~\kappa \frac{\partial u}{\partial n} \geq 0,~(u - g)\kappa \frac{\partial u}{\partial n} & = 0, \quad \text{on $\partial \Omega$},
	\end{aligned}
\end{equation}
is often referred to as the scalar Signorini's or Poisson--Signorini problem~\cite{gustafsson2019nitsche}.
We note that problem~\eqref{eq:signorini} is semi-coercive:
the bilinear form associated to the left-hand side controls only
the $H^1$-seminorm, and the inequality constraint
is not sufficient to restore full coercivity.
We assume that $f$ and $g$ are such that a unique solution exists.

The corresponding continuous saddle point formulation of \eqref{eq:signorini} is: find $(u, \lambda) \in H^1(\Omega) \times H^{-1/2}(\partial \Omega)$, $\lambda \geq 0$, such that
\begin{equation}
	\label{eq:saddlepointcont}
	\mathcal{L}(u, \lambda) = \inf_{v \in H^1(\Omega)} \sup_{\mu \in H^{-1/2}(\partial \Omega),\,\mu \geq 0} \mathcal{L}(v, \mu),
\end{equation}
where the Lagrangian $\mathcal{L}$ was defined in \eqref{eq:lagrangian}. The connection between the two formulations is again given by $\lambda = \kappa \frac{\partial u}{\partial n}$.
After residual stabilization, the finite element formulation reads: find $(u_h, \lambda_h) \in V_h \times Q_h$, $\lambda_h \geq 0$, such that
\begin{equation}
	\label{eq:saddlepointfem}
	\mathcal{L}_h(u_h, \lambda_h) = \inf_{v_h \in V_h} \sup_{\mu_h \in Q_h,\,\mu_h \geq 0} \mathcal{L}_h(v_h, \mu_h),
\end{equation}
where the stabilized Lagrangian $\mathcal{L}_h$ was defined in \eqref{eq:stabilizedlagrangian}.
The error analysis of formulation \eqref{eq:saddlepointfem} was given in \cite{gustafsson2019nitsche}
where it is concluded, among other things,
that the formulation
is stable for any conforming $Q_h$.

The corresponding variational inequality can be written as: find $\lambda_h \in Q_h$, $\lambda_h \geq 0$, so that
\begin{align*}
    	 \int_\Omega \kappa \nabla u_h \cdot \nabla v_h \dx - \int_{\partial \Omega} \lambda_h v_h \ds  + \int_{\partial \Omega} \gamma  \Big(\lambda_h - \kappa \frac{\partial u_h}{\partial n}\Big)\, \kappa \frac{\partial v_h}{\partial n} \ds &= \int_\Omega fv_h \dx, \\
\int_{\partial \Omega} \Big(\kappa \frac{\partial u_h}{\partial n} - \frac1\gamma(u_h - g) - \lambda_h\Big)( \mu_h - \lambda_h) &\leq 0,
\end{align*}
for every $(v_h, \mu_h) \in V_h \times Q_h$, $\mu_h \geq 0$.
The second inequality implies that $\lambda_h$ is the orthogonal $L^2$ projection of 
$\kappa \frac{\partial u_h}{\partial n} - \frac1\gamma(u_h - g)$ onto the set of nonnegative functions in $Q_h$ \cite[Theorem 5.2]{brezis2011functional}.
Since $Q_h$ is arbitrary,
we may choose $Q_h \subset L^2(\partial \Omega)$ such that the orthogonal projection is given explicitly by the maximum operator, i.e.,
\begin{equation}
    \label{eq:disclagmultineq}
    \lambda_h = \max(\kappa \tfrac{\partial u_h}{\partial n} - \tfrac1\gamma(u_h - g), 0) =: (\kappa \tfrac{\partial u_h}{\partial n} - \tfrac1\gamma(u_h - g))_+.
\end{equation}

In order to obtain the Nitsche method, we again substitute the discrete Lagrange multiplier $\lambda_h$ from \eqref{eq:disclagmultineq} to \eqref{eq:saddlepointfem} and rearrange the terms to obtain: find $u_h \in V_h$ which minimizes
\begin{equation}
	\label{eq:nitscheineq}
     J(u_h) + \int_{\partial \Omega} \frac{\gamma}{2} \Big(\kappa  \frac{\partial u_h}{\partial n} - \frac{1}{\gamma}(u_h - g)\Big)_+^2\ds -  \int_{\partial \Omega} \frac{\gamma}{2} \Big(\kappa \frac{\partial u_h}{\partial n}\Big)^2 \ds.
\end{equation}
Note that the corresponding method for the equality constraint
can be obtained by dropping the maximum operator (subscript $+$).
\section{Structure of the method}
\label{sec:nitschestructure}

We now generalize the method to arbitrary inequality constraints of the form $\beta(u_h) \geq 0$ where $\beta$ is an affine function.
We restrict our definition to affine functions because the existing a priori error estimates for similar methods, that we are aware of,
rely on the affinity of $\beta$; see, e.g., \cite{gustafsson2017mixed, gustafsson2017nitsche, gustafsson2019error}.
A discussion about extending the method
to nonlinear constraints is given at the end of
Section~\ref{sec:conclusions}.

The method \eqref{eq:nitscheineq} has the following main ingredients: the energy of the source problem, i.e., $J$; the constraint to impose, i.e., $\beta(u_h) = u_h-g \geq 0$; the corresponding definition of the continuous Lagrange multiplier, i.e., $\lambda(u_h) = \kappa \frac{\partial u_h}{\partial n}$; and the correct scaling of the stabilization parameter $\gamma$ with respect to the mesh parameter $h$ and any material parameters, i.e., $\gamma(h, \kappa)|_K=\tfrac{\alpha h_K}{\kappa}$.
In conclusion, the general form of the Nitsche method can be written as follows: find $u_h \in V_h$ that minimizes the functional
\begin{equation}
	\label{eq:nitschegen}
	J(u_h) + \int_{\Gamma} \frac{\gamma(h,\kappa)}{2} \Big(\lambda(u_h) - \frac{1}{\gamma(h,\kappa)}\beta(u_h)\Big)_+^2\ds -  \int_{\Gamma} \frac{\gamma(h, \kappa)}{2} \lambda(u_h)^2 \ds.
\end{equation}
The final method is formulated as a minimization problem instead of a variational one as it is the starting point for our implementation of the method; see Section~\ref{sec:novel} and \ref{app:newton}.
In particular, the calculation of the second functional derivative of \eqref{eq:nitschegen}
is performed automatically in our computational implementation of the Nitsche method.

\begin{remark}
	We  emphasize that the general form \eqref{eq:nitschegen} 
	was obtained
	following the derivation of \eqref{eq:nitscheineq}.
	Equivalently, it can be deduced starting
	from a similar minimization form existing in
	the literature  for other physical problems;
	see, e.g., \cite[p.~135]{chouly2013nitsche} for elastic
	contact problems and \cite[p.~364]{burman2017galerkin} for the obstacle problem.
	The form \eqref{eq:nitschegen} is also present in \cite{burman2023augmented}
	where the authors discuss an alternative approach using
	augmented Langrangian methods.
\end{remark}


Hence, for  new problems written in the general form \eqref{eq:nitschegen}, with  energy $J(u_h)$ and  finite element space $V_h$, it is necessary to define:
\begin{enumerate}
	\item[1.] A constraint $\beta = \beta(u_h) \geq 0$ and the subset $\Gamma \subset \Omega$ (or $\Gamma\subset \partial\Omega$) where the constraint
	      can become active.
	\item[2.] An expression for the continuous Lagrange multiplier (i.e., the contact force) as $\lambda = \lambda(u_h)$.
	\item[3.] An expression for the parameter $\gamma = \gamma(h, \kappa)$.
\end{enumerate}

In order to find an expression for $\lambda = \lambda(u_h)$, it is necessary to derive the strong formulation of the continuous saddle point problem.
This is usually obtained by first deriving the
variational formulation and then integrating by parts
to find the corresponding differential equations and boundary conditions,
some of which will depend directly on $\lambda$.
For contact-type constraints between two domains, 
prior research suggests that the parameter $\gamma$ and, hence, 
the expression for $\lambda$ should be defined on the 
less stiff side; cf., e.g.,~\cite{gustafsson2019error}.

Physically, equation \eqref{eq:nitschegen} shows that an important role of $\gamma$ is to scale the units of displacement in $\beta$ to the units of force in $\lambda$.
More precisely, $\lambda$ should be proportional to $\frac{\beta}{\gamma}$ with the proportionality constant independent of $h$ and $\kappa$.
Therefore,  the material parameters
of the problem must be included in $\gamma$.
Mathematically speaking,
a correct scaling of $\gamma$
is necessary to
prove the uniform stability of the
discrete saddle point formulation
using
an \emph{inverse inequality}; cf., e.g., \cite{gustafsson2017mixed, gustafsson2017nitsche, gustafsson2019error}.

Finally, we wish to remark that the penalty method
can be recovered from this formulation by removing 
the Lagrange multiplier terms (formally, $\lambda=0$), 
and that the Nitsche method for equality constraints
is obtained simply by dropping the maximum operator (subscript $+$).
Additionally, nonpositive constraints can be accommodated by changing 
the maximum operator to the minimum operator. 
The solution of the minimization problem \eqref{eq:nitschegen} using Newton's method
is briefly explained in Appendix~A with a link to our numerical implementations \cite{sourcepackage}.




\section{Some existing methods}

\subsection{Membrane obstacle problem}

Let $\Omega$ be a polygonal/polyhedral domain
representing a Poisson membrane in its
undeformed state.
Minimization of the energy
\begin{equation}
	J(u) = \frac12\int_\Omega \kappa \nabla u \cdot \nabla u \dx - \int_\Omega f u \dx
\end{equation}
corresponds to solving the Poisson equation $-\nabla \cdot (\kappa \nabla u) = f$ with $u \in H^1_0(\Omega)$.
Let now the deflection of the membrane $u$
be constrained by the rigid obstacle $g$
below the membrane, i.e., $u \geq g$.
The Lagrange multiplier corresponding
to the constraint $u - g \geq 0$ (in $\Omega$) is given by
$\lambda = -\nabla \cdot(\kappa \nabla u) - f =: -Lu - f$, cf.~\cite{gustafsson2017mixed, burman2017galerkin}.
The scaling is given by $\gamma|_K = \frac{\alpha h_K^2}{\kappa}$.
Hence, the Nitsche functional reads
\begin{equation}
	J(u_h) + \int_\Omega \frac{\gamma}{2}\left( - L_h u_h - f - \frac{1}{\gamma}(u_h - g) \right)_+^2 \dx - \int_\Omega \frac{\gamma}{2}\left( - L_h u_h- f\right)^2 \dx,
\end{equation}
where $L_h$ denotes an element-wise evaluation of $L$.

\subsection{Two-body linearized elastic contact}
Let $\Omega_1$ and $\Omega_2$ be two polygonal/polyhedral domains that are initially at
contact on $\Gamma = \partial \Omega_1 \cap \partial \Omega_2$.
The energy of the two-body linear elastic source problem reads
\begin{equation}
	\label{eq:2solid}
	J((\boldsymbol{w}_1, \boldsymbol{w}_2)) = \sum_{i=1}^2  \left\{  \frac12\int_{\Omega_i}\boldsymbol{\sigma}(\boldsymbol{w}_i) : \boldsymbol{\varepsilon}(\boldsymbol{w}_i)\dx -\int_{\Omega_i} \boldsymbol{f} \cdot \boldsymbol{w}_i \dx\right\},
\end{equation}
where
\begin{align}
	\label{eq:sigma}
	\bs{\sigma}(\bs{w}) & = 2 \mu_L\,\bs{\varepsilon}(\bs{w}) + \lambda_L \mathrm{tr}\,\bs{\varepsilon}(\bs{w}) \bs{I}, \\
	\label{eq:epsilon}
	\bs{\varepsilon}(\bs{w}) & = \frac12(\nabla \bs{w} + \nabla \bs{w}^T),
\end{align}
and $\bs{f}$ is a given body load
and $(\mu_L, \lambda_L)$ are the Lam\'e material parameters.

The physical non-penetration constraint on $\Gamma$
reads $\boldsymbol{u}_1 \cdot \boldsymbol{n}_1 + \boldsymbol{u}_2 \cdot \boldsymbol{n}_2 \leq 0$
where $\boldsymbol{n}_i$ is the outward normal vector on $\partial \Omega_i$, $i=1,2$.
The Lagrange multiplier corresponding
to the constraint $(\boldsymbol{u}_2 - \boldsymbol{u}_1) \cdot \boldsymbol{n}_1 \geq 0$
is given, e.g., by
$\lambda = -\boldsymbol{\sigma}(\boldsymbol{u}_1) \boldsymbol{n}_1 \cdot \boldsymbol{n}_1$.
Note that the choice to use the
normal traction on $\partial \Omega_1$ instead of $\partial \Omega_2$
is now arbitrary.  In general, the stabilization
should be on the less stiff body or the body with smaller elements; cf.~\cite{gustafsson2019error}.
The scaling is given by $\gamma|_K = \frac{\alpha h_K}{\mu_L}$.
The Nitsche functional \cite{chouly2015symmetric} reads
\begin{equation}
	\begin{aligned}
		J((\boldsymbol{u}_{1,h}, \boldsymbol{u}_{2,h})) & + \int_{\Gamma} \frac{\gamma}{2}\left( -\boldsymbol{\sigma}(\boldsymbol{u}_{1,h}) \boldsymbol{n}_1 \cdot \boldsymbol{n}_1- \frac{1}{\gamma}(\boldsymbol{u}_{2,h} \cdot \boldsymbol{n}_1 - \boldsymbol{u}_{1,h} \cdot \boldsymbol{n}_1) \right)_+^2 \dx \\
		& \quad - \int_{\Gamma} \frac{\gamma}{2}\left( -\boldsymbol{\sigma}(\boldsymbol{u}_{1,h}) \boldsymbol{n}_1 \cdot \boldsymbol{n}_1 \right)^2 \dx.
	\end{aligned}
\end{equation}



\section{Novel methods}
\label{sec:novel}

We now present some novel methods based on
formulation \eqref{eq:nitschegen}. More precisely, we consider the following problems:
\begin{enumerate}
	\item Two-membrane contact problem,
	\item Membrane against solid contact problem,
	\item Plate against plate contact problem,
	\item Kirchhoff plate with inequality boundary constraint.
\end{enumerate}
For each problem, we present the Nitsche formulation
and numerically investigate its convergence rate.
The numerical results are computed using scikit-fem~\cite{gustafsson2020scikit} which is used to solve the nonlinear variational problems using Newton's method.
For completeness, we describe how to manipulate the energy minimization problems for Newton's method in \ref{app:newton}. The programming implementation utilizes automatic differentiation with JAX \cite{bradbury2021jax} so that it is not necessary to manually calculate the functional derivatives involved in Newton's method. The source code is available in \cite{sourcepackage} for reproducing the numerical results.

Convergence rates are estimated via uniform mesh refinement.
In the absence of an analytical solution, we estimate
the convergence rate by evaluating the error between
two subsequent discrete solutions corresponding to mesh parameters
$2h$ and $h$.
Assuming that the error satisfies $\| u - u_h \| \leq C h^p$ in the energy norm for some $C>0$ and a given polynomial degree $p$ of the finite element basis,
we can show that
the difference between two subsequent discrete solutions should
be bounded by the same asymptotic rate:
$$\| u_h - u_{2h} \| \leq \| u - u_h \| + \| u - u_{2h} \| \leq C_1 h^p + C_2 (2h)^p \leq C_3 h^p, \quad C_1,C_2,C_3 > 0.$$

\subsection{Two-membrane contact problem}
\label{sec:2membrane}

\begin{figure}[H]
	\begin{center}
		\includegraphics[width=0.95\textwidth]{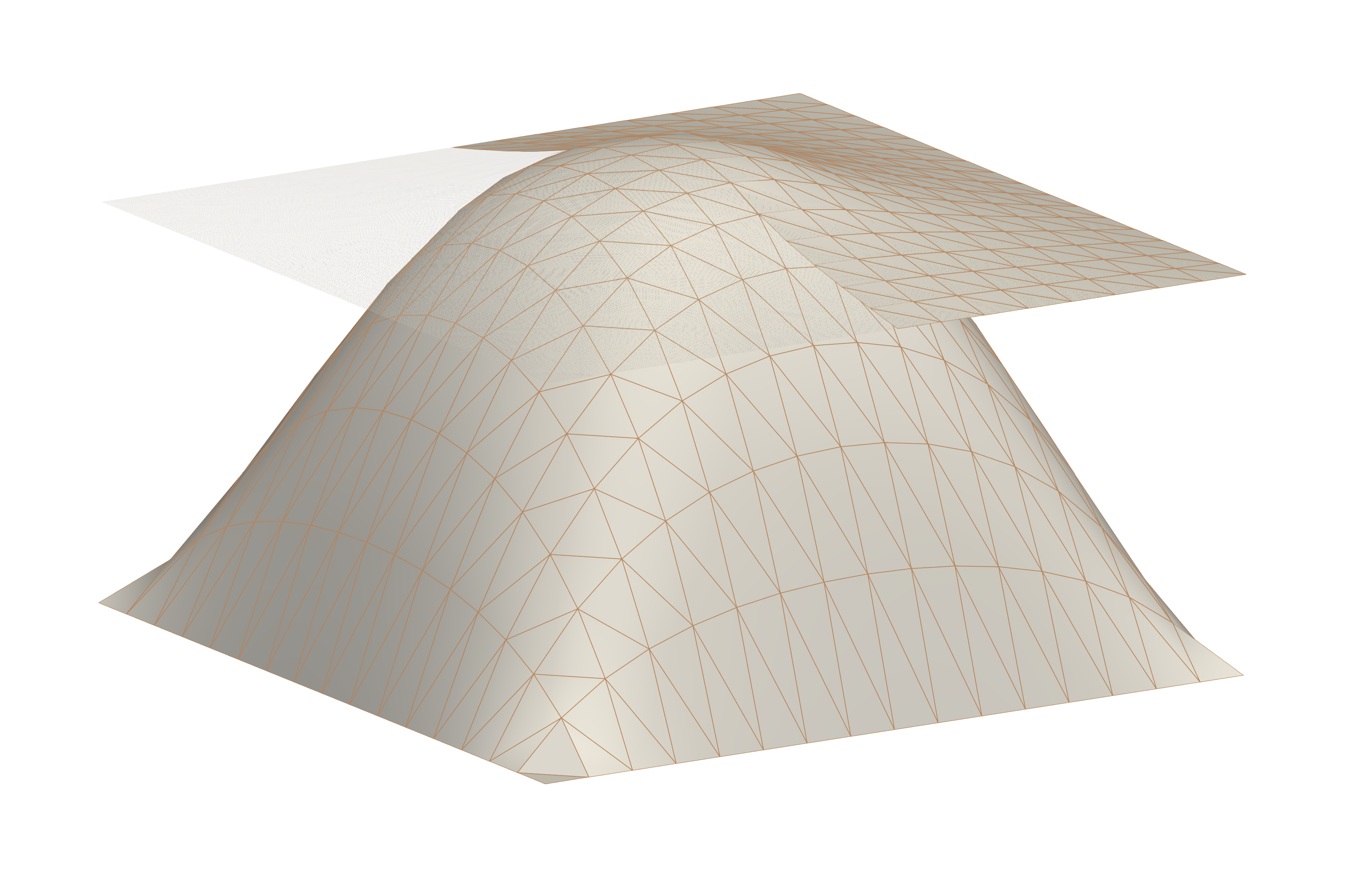}
	\end{center}
	\caption{Numerical solution for two membranes in contact using the Nitsche method with linear triangular elements.}\label{fig:2obstacle}
\end{figure}

Let $\Omega\subset \R^2$ be a polygonal domain describing the
shape of two structurally similar membranes 
modeled by the Poisson equation that are initially
$g$ (units of distance) apart. 
We denote the vertical displacements of the membranes by $u_1, u_2 : \Omega\to \R$ and their tensions by $\kappa_1, \kappa_2 > 0$, respectively. The membranes are initially separated by a gap $g>0$. 
The external loading $f_1, f_2 : \Omega\to \R$ can cause the membranes to
come in contact.
Therefore, we introduce a contact constraint $u_1-u_2\leq g$ in $\Omega$  and a corresponding contact pressure $\lambda: \Omega\to\R$ to prevent penetration.

The strong form of the problem reads as: find the displacements $u = (u_1, u_2)$  and the contact pressure $\lambda$ between the membranes such that 
\begin{equation}
	\label{eq:2membranestrong}
	\begin{alignedat}{2}
		- \nabla \cdot (\kappa_1 \nabla u_1) + \lambda & = f_1 && \qquad\text{in } \Omega,\\
		- \nabla \cdot (\kappa_2 \nabla u_2) - \lambda & = f_2 && \qquad\text{in } \Omega,\\
		u_1 = 0,\quad u_2 & = 0 && \qquad\text{on } \partial\Omega, \\
		u_1 - u_2 \leq g,\quad \lambda & \geq 0, \quad \lambda(u_2 - u_1 + g) = 0 && \qquad\text{in }\Omega. 
	\end{alignedat}
\end{equation}
\begin{remark}
	We define the problem such that the unknowns are the displacements of the membranes. Other authors, e.g. \cite{belgacem2009membranes, ben2012unilateral}, solve for the displaced coordinates of the membranes. 
\end{remark}
Explicitly, the second membrane is on top, initially $g$ units higher than the first membrane. The constraint $u_1-u_2\leq g$ prevents penetration by stating that the relative displacement of the membranes does not exceed the initial gap $g$. The contact pressure $\lambda$ appears with opposite signs in the equilibrium equations to represent the opposing forces in contact. The displacements $u_1$ and $u_2$ are zero on the boundary and we shall numerically approximate them in $\Omega$. We analytically solve and substitute the discrete $\lambda$ in the Nitsche method.

To derive the Nitsche method, we first identify the total potential energy of the unconstrained system, i.e.~ignoring the contact force $\lambda$. Multiplying the first equilibrium equation with a test function $v_1\in H_0^1(\Omega)$ and integrating by parts yields the weak form
\begin{align*}
	\int_\Omega -\nabla \cdot (\kappa_1\nabla u_1)v_1 = \int_\Omega \kappa_1\nabla u_1\cdot \nabla v_1 & = \int_\Omega f_1 v_1, \\ 
\implies \int_\Omega\kappa_1 \nabla u_1\cdot \nabla v_1 -\int_\Omega f_1 v_1& = 0.
\end{align*}
It is clear that the first term is the variation (Gateaux derivative) of $\frac{1}{2}\int_\Omega\kappa_1 \nabla u_1 \cdot \nabla u_1$. The same steps hold for the second equation. Hence, summing the two equations results in the total energy for the two-membrane source problem 
\begin{equation}
	J((u_1, u_2)) = \frac{1}{2}\int_\Omega \kappa_1\nabla u_1 \cdot \nabla u_1 \dx - \int_\Omega f_1 u_1 \dx + \frac{1}{2}\int_\Omega \kappa_2\nabla u_2 \cdot \nabla u_2 \dx - \int_\Omega f_2 u_2 \dx.
\end{equation}

Next, we analyze the continuous saddle point problem and find the corresponding Lagrange multiplier as we move towards the general Nitsche formulation \eqref{eq:nitschegen}. 
The constraint is given by $\beta((u_1,u_2)) = u_2-u_1+g\geq 0$. 
Hence, for the Lagrange multiplier $\lambda((u_1,u_2))$,
consider the Lagrangian 
\begin{equation}
	\mathcal{L}((u_1, u_2), \lambda) = J((u_1, u_2)) - \int_\Omega \lambda((u_1, u_2))\beta((u_1,u_2))\dx.
\end{equation}
Let $\kappa_1 < \kappa_2$ such that the first membrane is less stiff. We can find the strong form of the Lagrange multiplier by setting the Gateaux derivative with respect to $v=(v_1, 0)$ to zero, and finally applying integration by parts:
\begin{align*}
	\mathcal{L}(u, \lambda; v) & = \bigg(\frac12 \int_\Omega \kappa_1\nabla(u_1 + \epsilon v_1)\cdot \nabla (u_1 + \epsilon v_1) \dx - \int_\Omega f_1 (u_1 + \epsilon v_1) \dx \\ 
							   & \quad\quad + \frac12 \int_\Omega \kappa_2\nabla u_2\cdot\nabla u_2\dx - \int_\Omega f_2 u_2\dx - \int_\Omega \lambda (u_2 - (u_1 + \epsilon v_1) + g)\dx \bigg)\bigg|_{\epsilon=0}  \\ 
							   & = \bigg(\frac12 \int_\Omega \kappa_1(\nabla u_1\cdot \nabla u_1 + 2\epsilon \nabla u_1 \cdot \nabla v_1 + \epsilon^2 \nabla v_1\cdot\nabla v_1) \dx - \int_\Omega f_1 u_1 + \epsilon f_1 v_1\dx \\ 
							   & \quad\quad + \frac12 \int_\Omega \kappa_2\nabla u_2\cdot\nabla u_2\dx - \int_\Omega f_2 u_2\dx - \int_\Omega \lambda(u_2 - u_1 - \epsilon v_1 + g)\dx\bigg)\bigg|_{\epsilon=0}  \\ 
							   & = \int_\Omega \kappa_1\nabla u_1\cdot \nabla v_1 - f_1 v_1 + \lambda v_1 \dx\\ 
							   & = \int_\Omega - \kappa_1\Delta u_1 v_1 \dx - f_1 v_1 + \lambda v_1 \dx \\ 
							   & = 0 \\ 
		\implies \lambda & = \kappa_1\Delta u_1 + f_1. 
\end{align*}

We now have the constraint and an expression for the Lagrange multiplier, but must still ensure that the stabilized discrete problem has the correct units. As described in Section~\ref{sec:nitschestructure}, $\gamma$ scales the discrete problem such that $\lambda \propto \frac{\beta}{\gamma}$. Because the discrete solution $u_h$ is differentiated twice in the definition $\lambda(u_h)$ and $\beta(u_h)$ includes the solution $u_h$ in its original units (i.e., in meters), $\gamma$ must here include a square of the mesh parameter. Further, as we derived $\lambda$ with respect to the less stiff side per \cite{gustafsson2019error}, the tension of the first membrane is included in $\lambda$ and correspondingly $\gamma$ must have its reciprocal. Hence, it must be that $\gamma(h,\kappa) = \frac{\alpha h^2}{\kappa_1}$ for the units to coincide, where $\alpha$ is dimensionless. Table \ref{tab:gamma_2membrane} summarizes the stabilization parameter scaling.

\begin{table}[H]
	\centering
    	\caption{Scaling of the stabilization parameter $\gamma$ for the two-membrane contact problem.}
	\begin{tabular}{lp{0.8\textwidth}}
		\toprule
		$\gamma(h,\kappa)$ & $\dfrac{\alpha h^2}{\kappa_1}$ \\[1em]
		Power of $h$      & $h^2$ \\[0.5em]
		Reasoning     & The multiplier $\lambda(u_h) = \kappa_1 \Delta u_{h,1} + f_1$ contains two derivatives of the solution, so $\gamma$ must include $h^2$ to balance the units of $\beta(u_h) = u_{h,2} - u_{h,1} + g$ in $\lambda(u_h) - \beta(u_h)/\gamma(h,\kappa)$, see \eqref{eq:nitschegen}. The scaling $1/\kappa_1$ appears because the stabilization is taken on the less stiff side, cf.~\cite{gustafsson2019error}. The coefficient $\alpha>0$ is dimensionless. \\
		\bottomrule
	\end{tabular}
	\label{tab:gamma_2membrane}
\end{table}

To collect, the Nitsche terms are 
\begin{align*}
	\lambda(u_h) = \kappa_1\Delta_h u_{h,1} + f_1, \quad 
	\beta(u_h) = u_{h,2} - u_{h,1} + g,  \quad
	\gamma(hz, \kappa) = \frac{\alpha h^2}{\kappa_1}, 
\end{align*}
where $\Delta_h$ denotes the element-wise Laplacian operator. 
Substituting these to \eqref{eq:nitschegen} produces the Nitsche problem: find $(u_{h,1}, u_{h,2})\in V_h\times V_h$ which minimizes
\begin{equation}
	\begin{aligned}
    \label{eq:2memenergy}
		J((u_{h,1}, u_{h,2})) & + \int_\Omega \frac{\alpha h^2}{2\kappa_1}\left( \kappa_1\Delta_h u_{h,1} + f_1 - \frac{\kappa_1}{\alpha h^2}(u_{h,2} - u_{h,1} + g) \right)_+^2 \dx \\ & - \int_\Omega \frac{\alpha h^2}{2\kappa_1}\left( \kappa_1\Delta_h u_{h,1}+ f_1\right)^2 \dx.
	\end{aligned}
\end{equation}
\begin{remark}
	The Nitsche terms and formulation are derived similarly for other problems in Section \ref{sec:novel}.
\end{remark}

\begin{figure}[H]
	\begin{center}
		\includegraphics[width=0.95\textwidth]{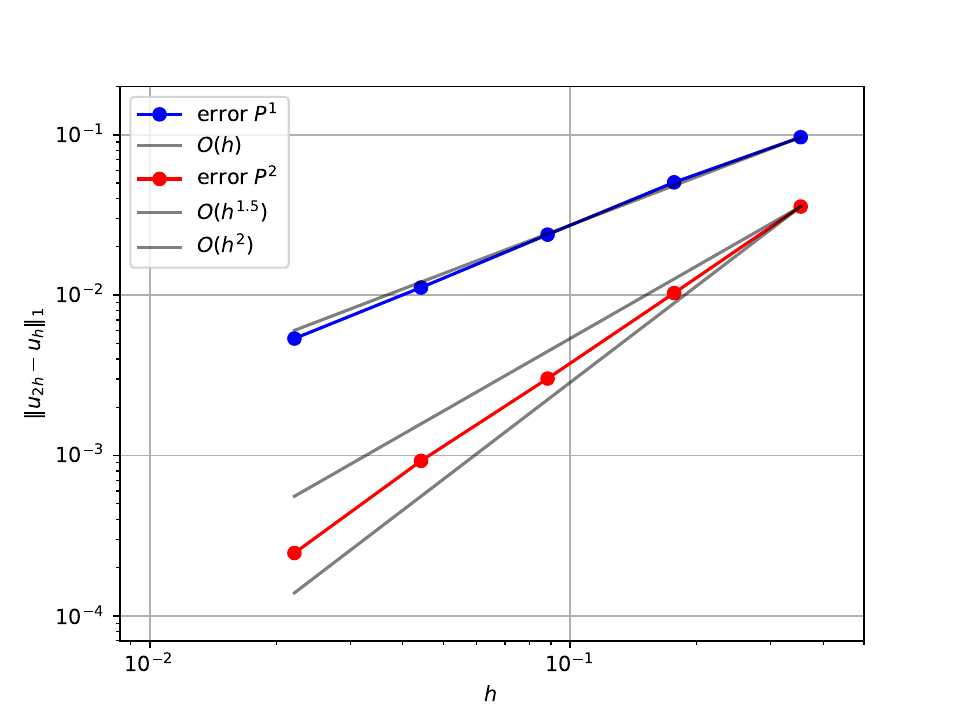}
	\end{center}
	\caption{Two membrane contact problem convergence rates follow the theoretical linear in the $H^1$ norm for linear element basis. 
	For quadratic basis, the convergence is slightly slower than quadratic. 
    Solutions to obstacle problems
    are known to be limited by
    $C^{1.1}$ regularity \cite{caffarelli1998obstacle}.
    This is because the second derivative, in general, has a jump discontinuity at the free boundary and, therefore, $u_1 \in H^{s}(\Omega)$, $s \leq 5/2$.
    Hence, standard a priori error analysis implies
    that the error is $O(h^{1.5})$.
    The estimated rate is slightly higher than $O(h^{1.5})$ which can be attributed to the matching nodes in the meshes of $u_1$ and $u_2$.}
    
	\label{fig:2obstacleconvergence}
\end{figure}

The problem was solved numerically over $\Omega=(0,1)^2$ with the parameterization $g=0.05$, $f_1=1, f_2=0, \kappa_1 = \kappa_2 = 1$ and 
$\alpha=10^{-2}$ 
using both linear and quadratic triangular elements. Note that the element-wise Laplacian is equal to zero for linear elements and constant for quadratic elements, which simplifies the formulation. The linear element numerical solution is displayed in Figure \ref{fig:2obstacle} and the approximated convergence rates in Figure \ref{fig:2obstacleconvergence}.


\begin{figure}
	\begin{center}
		\includegraphics[trim={2.1cm 1cm 2.5cm 1cm}, clip,width=\textwidth]{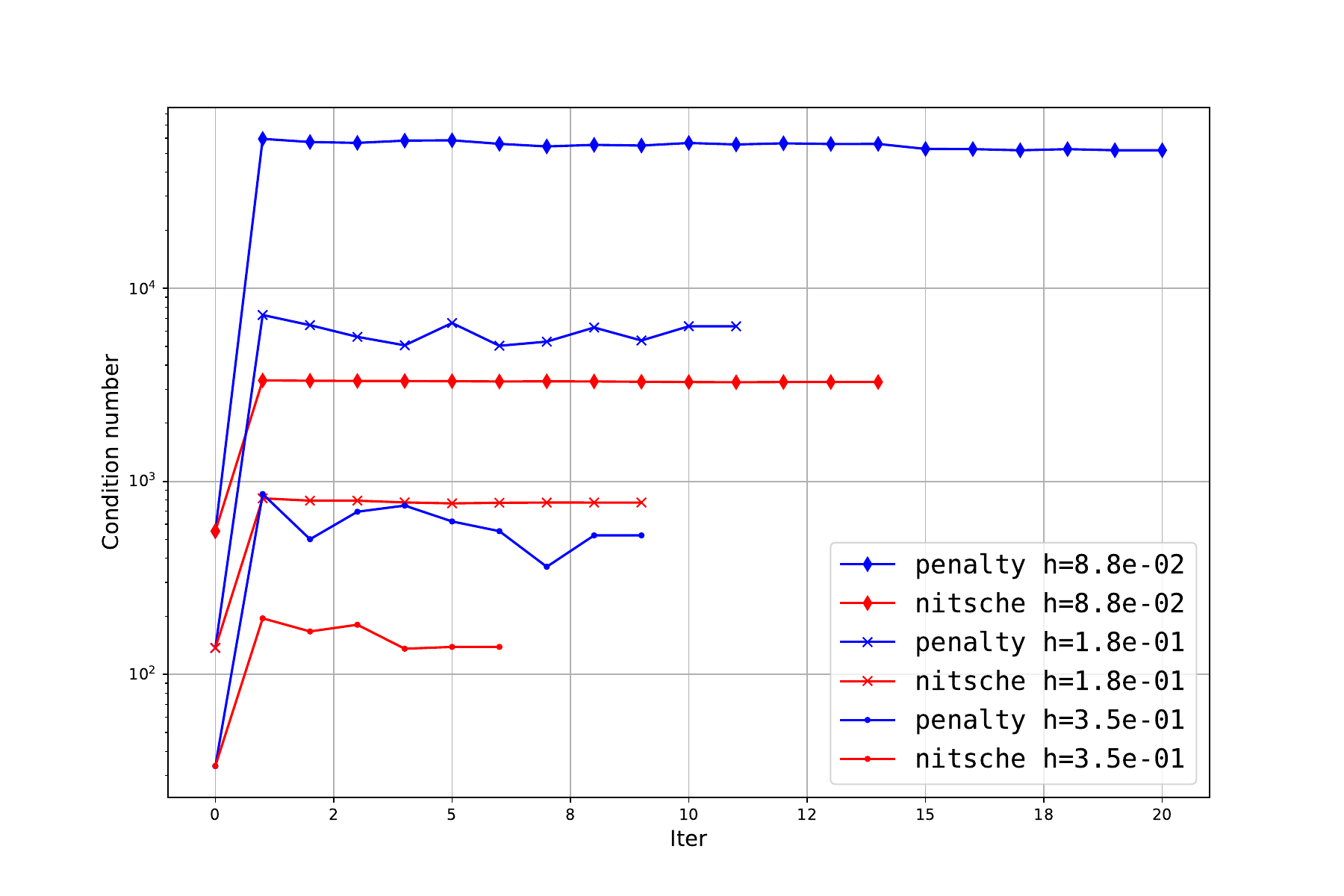}
	\end{center}
	\caption{Condition numbers of the Jacobian in the Newton iterations for the penalty and Nitsche variants of the two membrane contact problem with quadratic elements. The blue lines depict condition numbers of the penalty method and red lines of the Nitsche method, respectively. The markers distinguish problems of different sizes. Compared to the Nitsche method's $\gamma = \alpha h^2$, the penalty method requires $1/\epsilon = \alpha h^3$ to retain optimal convergence \cite{gustafsson2017finite} and has thus consistently larger conditioning numbers with slower convergence.}\label{fig:2obstacle_condnum}
\end{figure}

As described earlier, we can revert to the penalty method by removing the terms depending on $\lambda((u_1, u_2))$. However, the penalty method retains the optimal convergence rate only with linear elements, and to produce the same level of accuracy as Nitsche for quadratic or higher order elements we need to increase the power of $h$ in $1/\epsilon$, cf.~\cite{gustafsson2017finite}. This leads to higher condition numbers in the Newton iteration. Figure \ref{fig:2obstacle_condnum} provides an example in the quadratic case.

\subsection{Membrane against elastic solid}
\label{sec:membranesolid}

Let $\Omega_1$ describe the membrane
$$\Omega_1 = \{ (x, y, 0) \in \mathbb{R}^3 : 0 < x, y < 1 \}$$
and
$\Omega_2$ describe the linear elastic cube, initially $g$ units apart from $\Omega_1$:
$$\Omega_2 = \{ (x, y, z + g) \in \mathbb{R}^3 : 0 < x, y, z < 1 \}.$$
The energy of the source problem including membrane  and linear elastic
solid reads
\begin{equation}
	J((u_1, \boldsymbol{u}_2)) = J_1(u_1) + J_2(\boldsymbol{u}_2)
\end{equation}
where
\begin{equation}
	J_1(u_1) = \int_{\Omega_1} \nabla u_1 \cdot \nabla u_1 \dx - \int_{\Omega_1} f u_1 \dx
\end{equation}
and
\begin{equation}
	J_2(\boldsymbol{u}_2) =   \int_{\Omega_2}\boldsymbol{\sigma}(\boldsymbol{u}_2) : \boldsymbol{\varepsilon}(\boldsymbol{u}_2)\dx.
\end{equation}
The functions $\bs{\sigma}$ and $\bs{\varepsilon}$ have been defined in equations $\eqref{eq:sigma}$ and $\eqref{eq:epsilon}$, respectively.

Given the initial gap $g>0$, we set the physical non-interpenetration constraint as $u_1 - \bm{u}_2\cdot \bm n \leq g$ (in $\Omega_1$) which results in $\beta((u_1, \bm{u}_2)) = \bm u_2 \cdot \bm n - u_1 + g \geq 0$. The Lagrange multiplier for the problem is 
$\lambda((u_1, \bm u_2)) = - \bm \sigma(\bm u_2)\bm n_2 \cdot \bm n_2$
and the scaling is $\gamma(h) = \alpha h$. Table~\ref{tab:gamma_membranesolid} summarizes the stabilization parameter scaling.

\begin{table}[H]
	\centering
    \caption{Scaling of the stabilization parameter $\gamma$ for the membrane against elastic solid problem.}
	\begin{tabular}{lp{0.8\textwidth}}
		\toprule
		$\gamma(h)$    & $\alpha h$ \\[0.5em]
		Power of $h$   & $h$ \\[0.5em]
		Reasoning  & The multiplier $\lambda(u_h) = -\bm\sigma(\bm u_{h,2})\bm n_2 \cdot \bm n_2$ contains one derivative of the solution, so $\gamma$ must include $h$ to balance the units of $\beta(u_h) = \bm u_{h,2}\cdot\bm n - u_{h,1} + g$ in $\lambda(u_h) - \beta(u_h)/\gamma(h)$, see \eqref{eq:nitschegen}. The elastic moduli are normalized to unity in this example. The coefficient $\alpha>0$ is dimensionless. \\
		\bottomrule
	\end{tabular}
	\label{tab:gamma_membranesolid}
\end{table}

\begin{figure}[H]
	\begin{center}
		\includegraphics[width=0.95\textwidth]{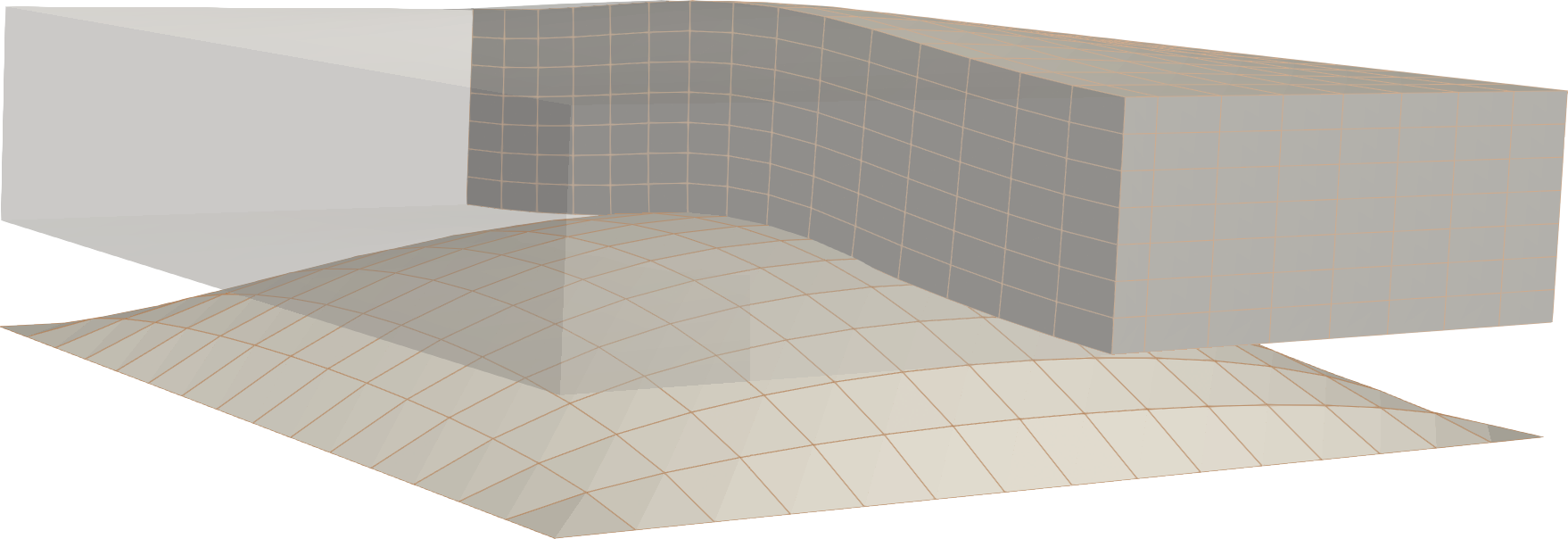}
	\end{center}
	\caption{Numerical solution for the membrane against elastic solid problem using the Nitsche method with linear hexahedral elements for the solid and linear quadrilateral elements for the membrane.}\label{fig:membranesolid}
\end{figure}

\begin{figure}[H]
	\begin{center}
		\includegraphics[width=0.95\textwidth]{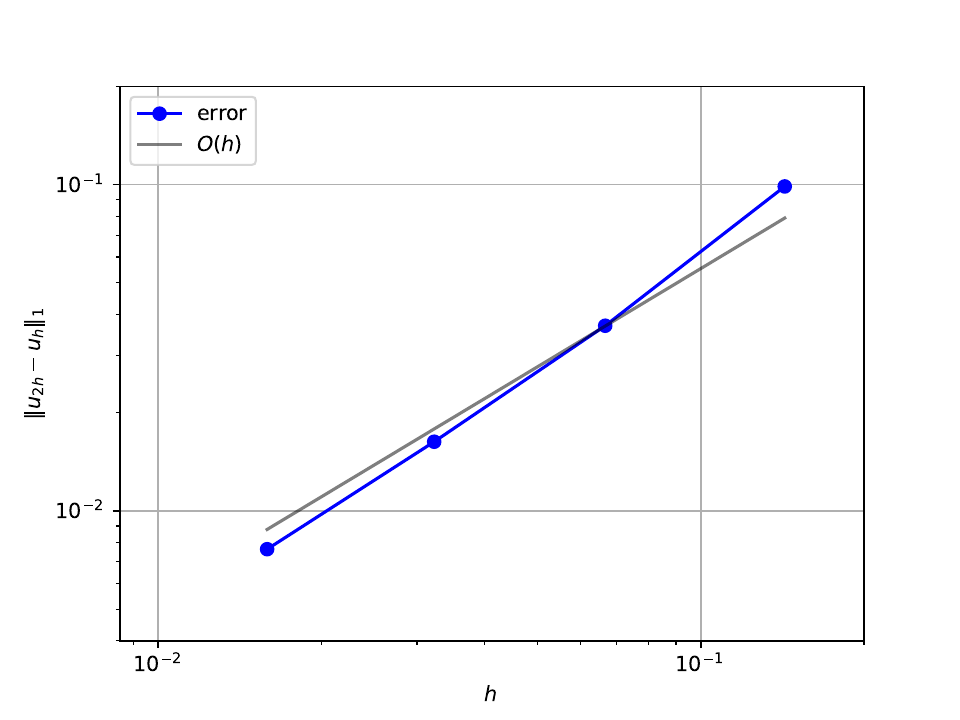}
	\end{center}
	\caption{ The convergence rate of the membrane against elastic solid problem follows the theoretical linear convergence in the $H^1$ norm with linear elements. }
	\label{fig:membranesolidconvergence}
\end{figure}

\begin{figure}
	\centering
	\begin{subfigure}[b]{0.49\textwidth}
		\centering
		\includegraphics[width=\textwidth]{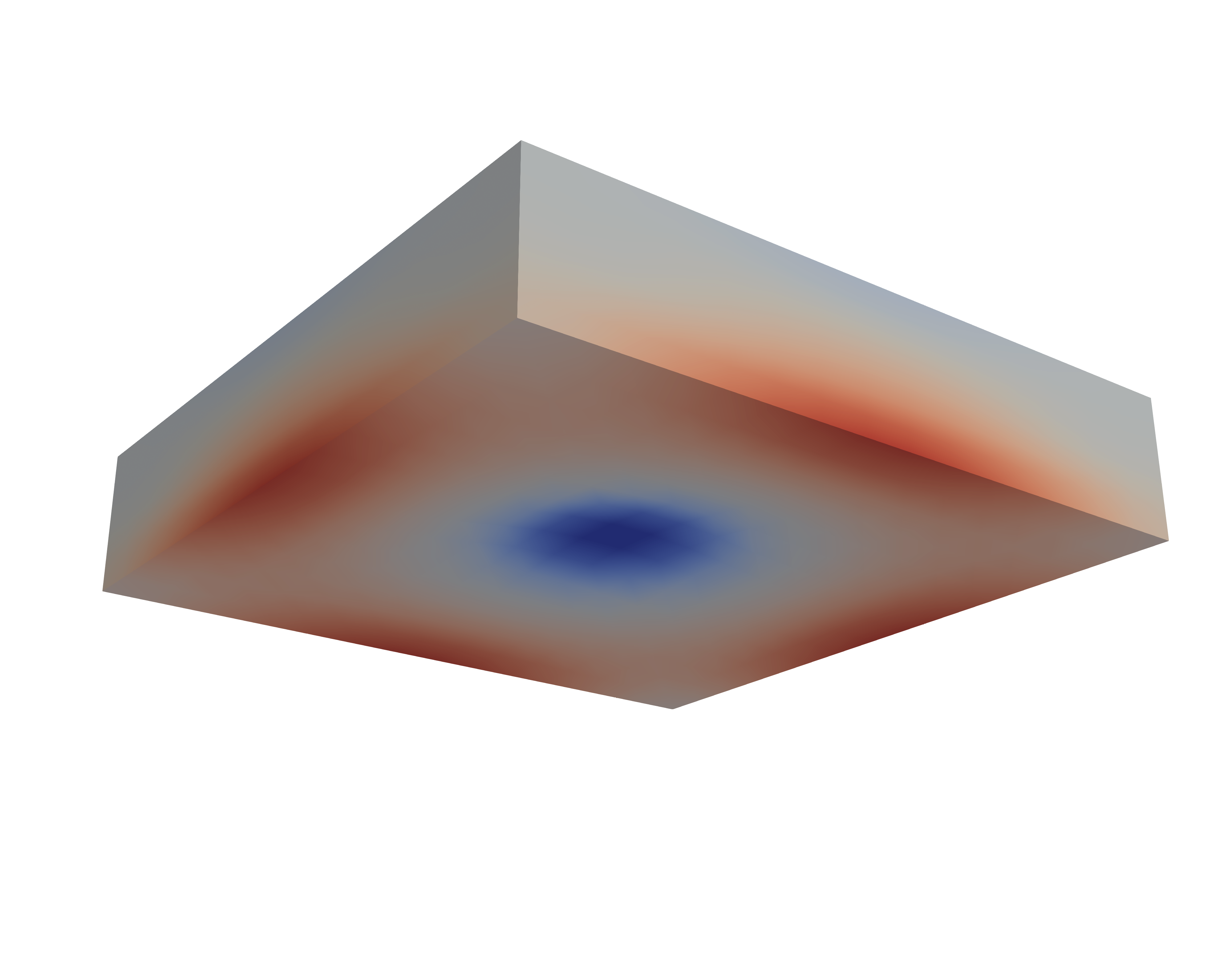}
		\caption[]%
		{{\small }}
	\end{subfigure}
	\hfill
	\begin{subfigure}[b]{0.49\textwidth}
		\centering
		\includegraphics[width=\textwidth]{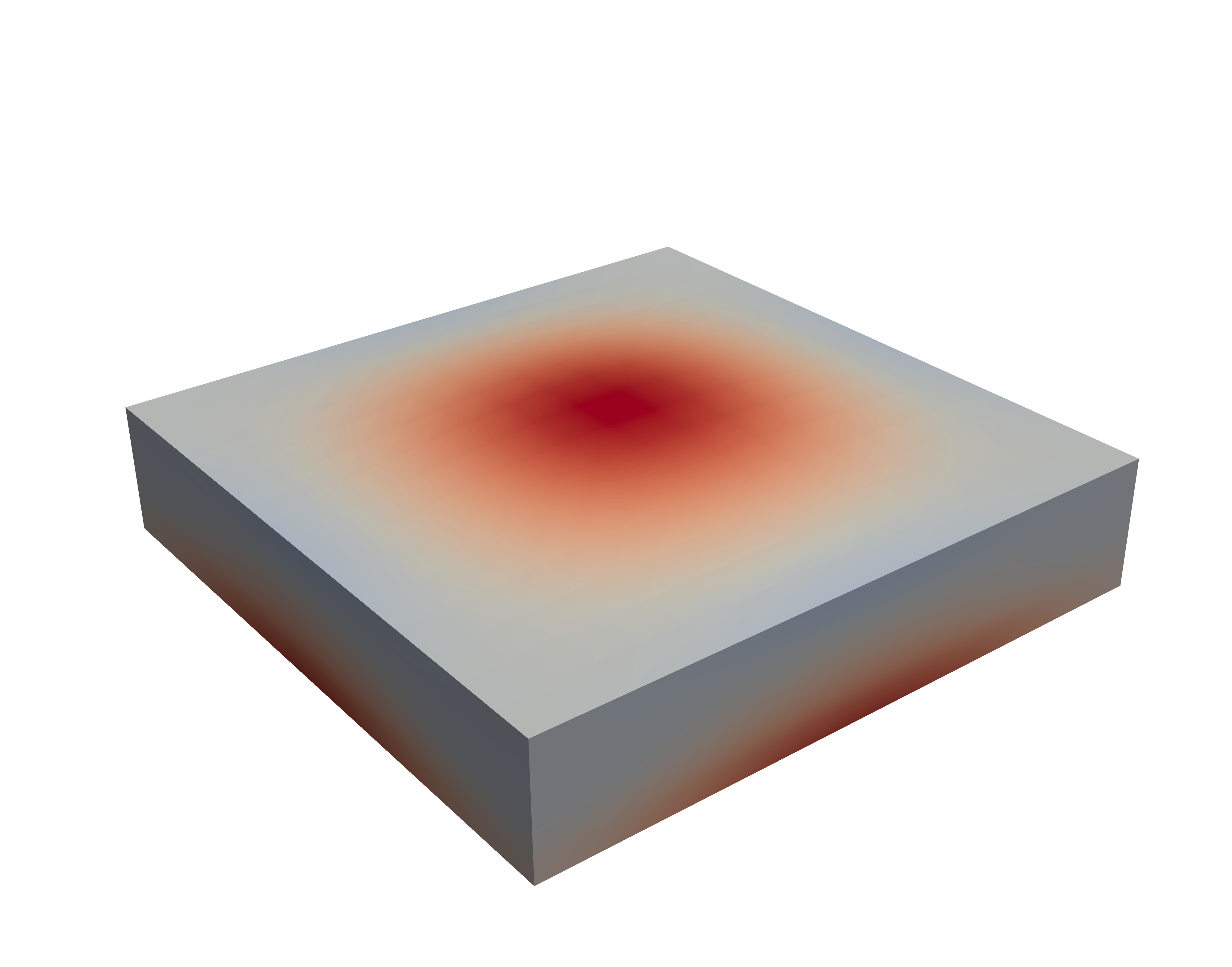}
		\caption[]%
		{{\small }}
	\end{subfigure}
	\caption[ Lol ]
	{\small The first principal stress for the elastic solid in membrane against elastic solid problem.}
	\label{fig:membranesolid_principal_stress}
\end{figure}

Thus, the Nitsche method minimizes the constrained energy
\begin{equation}
	\begin{aligned}
	J(u_{1,h}, \boldsymbol{u}_{2,h}) & + \int_{\Omega_1} \frac{\alpha h}{2} \left( - \bm \sigma(\bm u_2)\bm n_2 \cdot \bm n_2 - \frac{1}{\alpha h} (\bm u_2 \cdot \bm n - u_1 + g)\right)_+^2 \\ 
									 & - \int_{\Omega_1} \frac{\alpha h}{2} (- \bm \sigma(\bm u_2)\bm n_2 \cdot \bm n_2)^2 \dx.
	\end{aligned}
\end{equation}
The problem was solved numerically using a gap of $g=0.1$, $\mu_L = \lambda_L=1$, $f=2$ and 
$\alpha=10^{-2}$ 
using linear hexahedral elements for the solid and linear quadrilateral elements for the membrane. The numerical solution is displayed in Figure \ref{fig:membranesolid} and the approximated linear convergence rate in Figure \ref{fig:membranesolidconvergence}. The first principal stress for the solid is presented in Figure \ref{fig:membranesolid_principal_stress}.


\vspace{0.5cm}

\subsection{Plate against plate}
\label{sec:2plate}


The plate against plate problem is a Kirchhoff--Love plate model with clamped plates, which 
follows very closely the two-membrane contact problem with the difference of higher order derivatives in the strain energy. The energy of the source problem is
\begin{equation}
	J((u_1, u_2)) = \frac{1}{2}\int_\Omega \sum_{i,j=1}^2 \left(\frac{\partial u_1}{\partial x_i x_j}\right)^2 \dx - \int_\Omega f_1 u_1 \dx + \frac{1}{2}\int_\Omega \sum_{i,j=1}^2 \left(\frac{\partial u_2}{\partial x_i x_j}\right)^2 \dx - \int_\Omega f_2 u_2 \dx.
\end{equation}
We assume again a gap $g>0$ such that $u_1 - u_2 \leq g$ in $\Omega$ for the constraint $\beta{(u_1,u_2)} = u_2 - u_1 + g \geq 0$. Then, we have $\lambda((u_1, u_2))=-\Delta^2 u_1 - f_1$ and $\gamma(h) = \alpha h^4$. Table~\ref{tab:gamma_plateplate} summarizes the stabilization parameter scaling.

\begin{table}[H]
	\centering
    	\caption{Scaling of the stabilization parameter $\gamma$ for the plate against plate problem.}
	\begin{tabular}{lp{0.8\textwidth}}
		\toprule
		$\gamma(h)$    & $\alpha h^4$ \\[0.5em]
		Power of $h$   & $h^4$ \\[0.5em]
		Reasoning  & The multiplier $\lambda(u_h) = -\Delta^2 u_{h,1} - f_1$ contains four derivatives of the solution via the biharmonic operator, so $\gamma$ must include $h^4$ to balance the units of $\beta(u_h) = u_{h,2} - u_{h,1} + g$ in $\lambda(u_h) - \beta(u_h)/\gamma(h)$, see \eqref{eq:nitschegen}. The plate bending moduli are normalized to unity in this example. The coefficient $\alpha>0$ is dimensionless. \\
		\bottomrule
	\end{tabular}
	\label{tab:gamma_plateplate}
\end{table}

Hence, the Nitsche problem minimizes the energy
\begin{equation}
	\begin{aligned}
		J((u_{h,1}, u_{h,2})) & + \int_\Omega \frac{\alpha h^4}{2}\left( - \Delta_h^2 u_{h,1} - f_1 - \frac{1}{\alpha h^4}(u_{h,2} - u_{h,1} + g) \right)_+^2 \dx \\ & - \int_\Omega \frac{\alpha h^4}{2}\left( - \Delta_h^2 u_{h,1}- f_1\right)^2 \dx,
	\end{aligned}
\end{equation}
where $\Delta_h^2$ denotes the element-wise biharmonic operator. Note that the Nitsche method for plate obstacle problems were considered in~\cite{fabre2021nitsche}, which suggests that the Kirchhoff model can be suitable as long as 
the transverse compression of the plate can be neglected.

\begin{figure}[H]
	\begin{center}
		\includegraphics[width=0.85\textwidth]{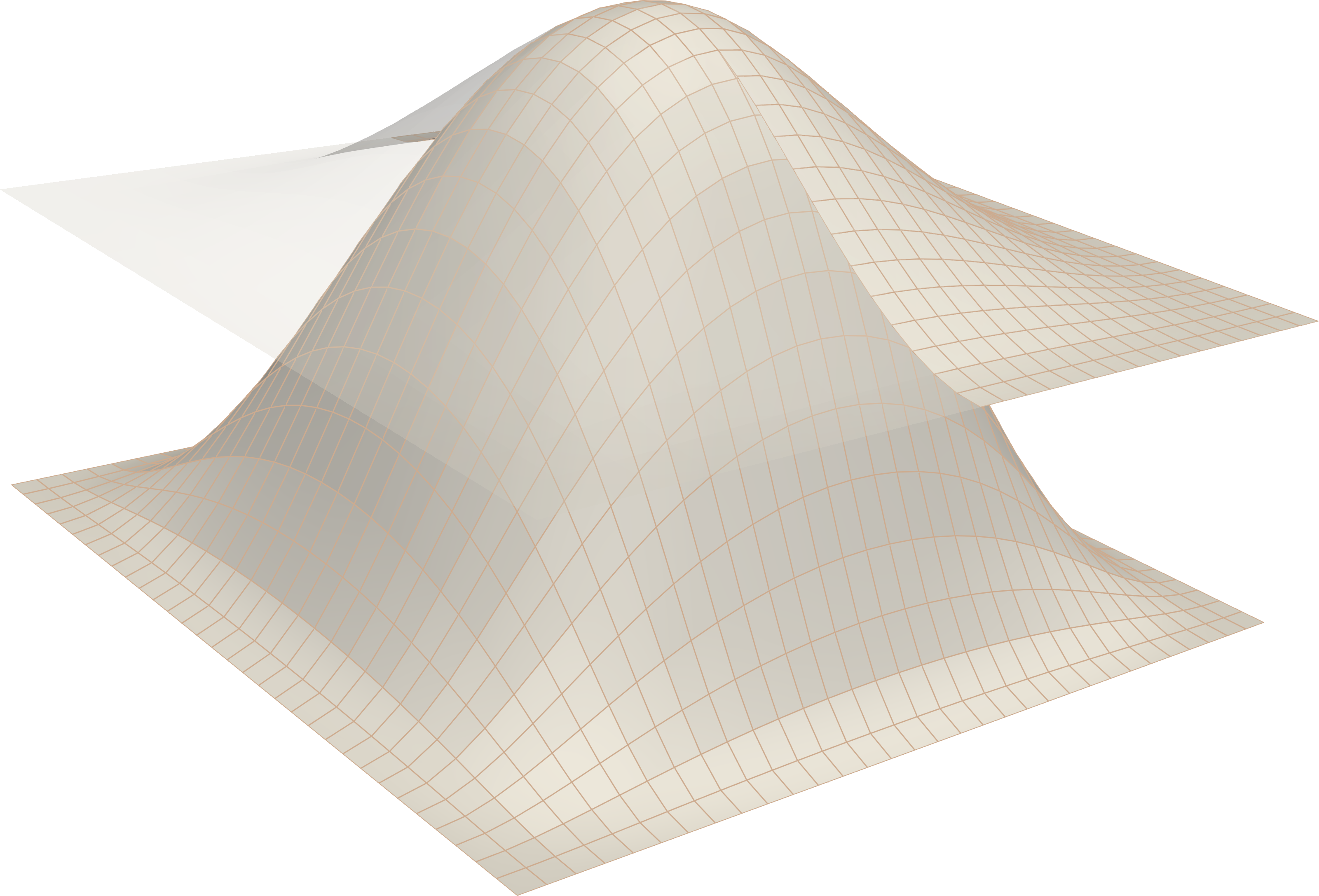}
	\end{center}
	\caption{Numerical solution for two plates in contact using the Nitsche method with quadrilateral BFS elements. 
	}\label{fig:2plate}
\end{figure}

\begin{figure}[H]
	\begin{center}
		\includegraphics[width=0.8\textwidth]{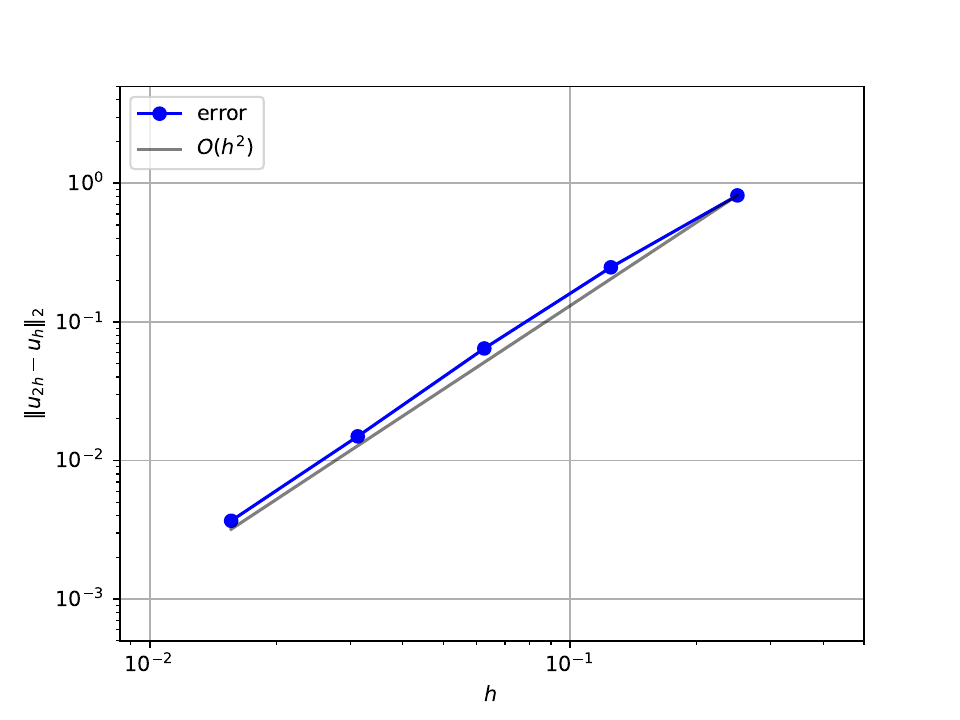}
	\end{center}
	\caption{Two plate contact problem convergence rate follows the theoretical quadratic convergence in the $H^2$ norm with quadrilateral BFS elements \cite{bogner1965generation}.}\label{fig:2plateconvergence}
\end{figure}

The problem was solved numerically over $\Omega=(0,1)^2$, a gap of $g=0.05$, $f_1=100$, $f_2=0$ and 
$\alpha=10^{-2}$
using quadrilateral Bogner--Fox--Schmit elements~\cite{bogner1965generation}.
The numerical solution is displayed in Figure \ref{fig:2plate} and the approximated quadratic convergence rate in Figure \ref{fig:2plateconvergence}.

\subsection{Kirchhoff plate with inequality boundary condition}
\label{sec:platecorner}

Let $\Omega = (0,1)^2$.
For the Kirchhoff plate with a simply-supported inequality boundary condition~\cite{nazarov2012hinged}, the energy with normalized units reads
\begin{equation}
	J(u) = \int_\Omega \sum_{i,j=1}^2 \left(\frac{\partial^2 u}{\partial x_i\partial x_j}\right)^2 - \int_\Omega fu.
\end{equation}
We consider the problem where $\Gamma = \partial\Omega$ and $\beta(u) = u \geq 0$ on $\Gamma$. Then, the Lagrange multiplier is $\lambda(u) = V_n(u) := Q_n + \frac{\partial M_{ns}}{\partial s}$, where $V_n(u)$ is called the \emph{Kirchhoff shear force} and
$$M(u) = \frac{Ed^3}{12(1+\nu)}\left(K(u) + \frac{\nu}{1-\nu}(\text{tr}\,K(u))I\right),$$ 
where $K(u) = -\nabla \nabla u$ is the curvature, $d$ denotes the plate thickness, $I$ is the identity tensor, and $E$ and $\nu$ are Young’s modulus and Poisson ratio, respectively. Further,
$Q= \text{div}\: M$, $n$ is the normal vector and $s$ the tangent vector on $\Gamma$. For more details, we refer to~\cite{gustafsson2021nitsche}. As mentioned, we assume normalized units so that $\nu=0$ and $\frac{Ed^3}{12}=1$ and it follows that the scaling parameter is $\gamma(h) = \alpha h^3$. Table~\ref{tab:gamma_platecorner} summarizes the stabilization parameter scaling.

\begin{table}[H]
	\centering
    	\caption{Scaling of the stabilization parameter $\gamma$ for the Kirchhoff plate with inequality boundary condition.}
	\begin{tabular}{lp{0.8\textwidth}}
		\toprule
		$\gamma(h)$    & $\alpha h^3$ \\[0.5em]
		Power of $h$   & $h^3$ \\[0.5em]
        Reasoning  & The multiplier $\lambda(u_h) = V_n(u_h) = Q_n + \partial M_{ns}/\partial s$ contains three derivatives of the solution, so $\gamma$ must include $h^3$ to balance the units of $\beta(u_h) = u_h$ in $\lambda(u_h) - \beta(u_h)/\gamma(h)$, see \eqref{eq:nitschegen}. The plate moduli are normalized to unity in this example. The coefficient $\alpha>0$ is dimensionless.  \\
		\bottomrule
	\end{tabular}
	\label{tab:gamma_platecorner}
\end{table}

The Nitsche functional is thus
\begin{align}
	J(u_h) + \int_{\Gamma} \frac{\alpha h^3}{2} \Big(V_{n,h}(u) - \frac{1}{\alpha h^3}u_h\Big)_+^2\ds -  \int_{\Gamma} \frac{\alpha h^3}{2} V_{n,h}(u)^2 \ds,
\end{align}
where $V_{n,h}$ is the element-wise evaluation of $V_n$.

\begin{figure}[H]
	\begin{center}
		\includegraphics[trim={2.5cm 2cm 1cm 2cm}, clip,width=\textwidth]{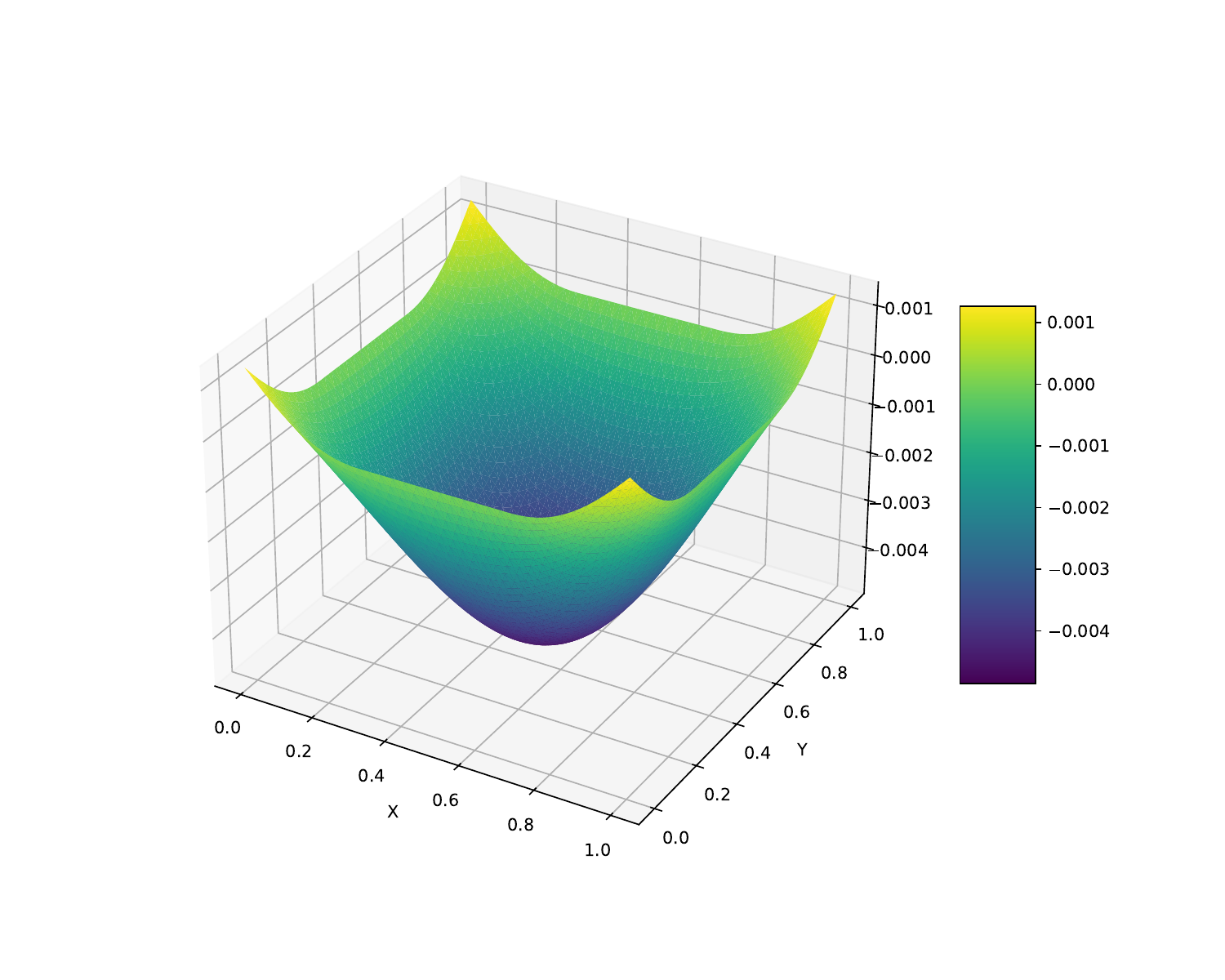}
	\end{center}
	\caption{The corners of the plate are displaced upwards when downwards point load is applied at the middle of the plate, approximate solution using Nitsche's method with quadrilateral BFS elements.}\label{fig:platecorner}
\end{figure}

\begin{figure}[H]
	\begin{center}
		\includegraphics[width=0.95\textwidth]{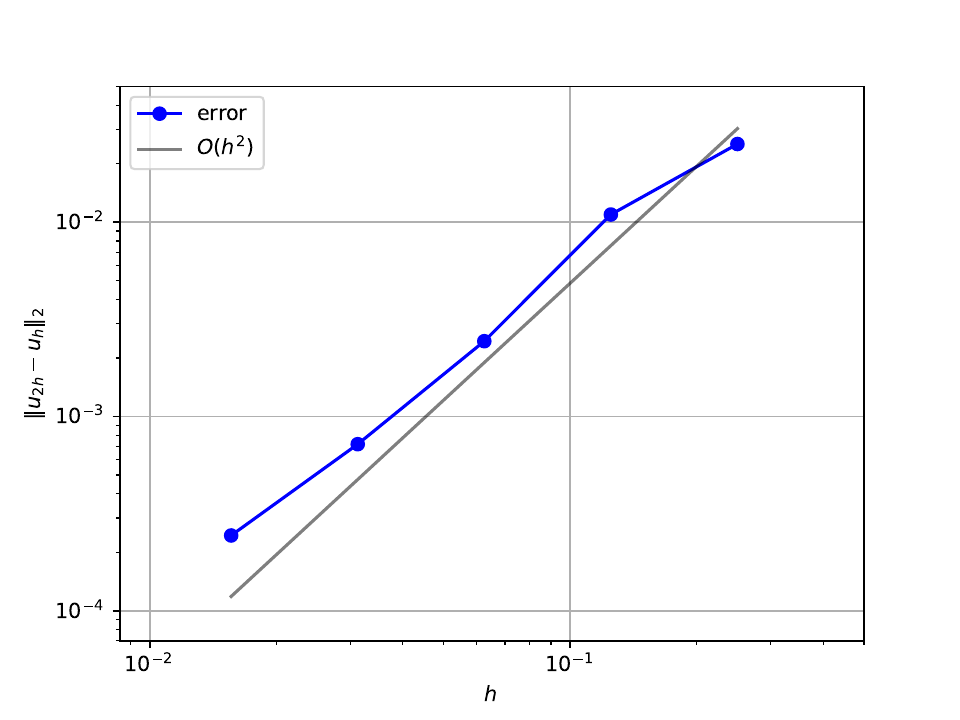}
	\end{center}
	\caption{The convergence rate of the Kirchhoff plate with inequality boundary condition follows the theoretical quadratic convergence in the $H^2$ norm with quadrilateral BFS elements \cite{bogner1965generation}.}\label{fig:platecornerconvergence}
\end{figure}

The problem was solved numerically with $f=-1$ and $\alpha=1/4$ using quadrilateral Bogner--Fox--Schmit elements \cite{bogner1965generation}. The numerical solution is displayed in Figure \ref{fig:platecorner} and the approximated quadratic convergence rate in Figure \ref{fig:platecornerconvergence}.

\section{Conclusions}

\label{sec:conclusions}

While the Nitsche method is commonly viewed as a consistency correction to the penalty method, we demonstrate how new methods can be obtained starting from a general minimization form. This allows for developing novel Nitsche methods for problems with both affine equality and inequality constraints in a simple manner. We derive and implement several new problems using the approach. Our numerical experiments suggest that the methods converge optimally while the more careful numerical analysis of the general formulation remains a topic of future work.

\section*{Limitations}

We have restricted the constraint $\beta$ to affine functions because the
existing error analyses for specific implementations of the method rely on the affinity of $\beta$.
It is clear that formally the Nitsche method
can be defined for any, even nonlinear, constraint $\beta$.
While we are not aware of
any results or attempts on the error analysis of such methods, we have included
exploratory numerical evidence with a nonlinear constraint in Appendix~B.
For an example of Nitsche's method with a nonlinear constraint in the existing literature see, e.g., \cite{MLIKA2017265}.
More careful study of the presented general formulation with nonlinear
constraints is another topic of future work.

\section*{Acknowledgements}

We wish to thank Rolf Stenberg and the anonymous referees for valuable comments on the manuscript.

\appendix
\section{Numerical implementation}
\label{app:newton}
The discrete minimization problem~\eqref{eq:nitschegen} is solved
iteratively by Newton's method. Given an iterate $u_{k,h}\in U_h$, the
energy $J$ is linearized by computing its first and second Gateaux
derivatives, yielding the Newton step
$J''(u_{k,h};v_h,w_h)=-J'(u_{k,h};v_h)$ for all $v_h\in V_h$, which is linear
in the increment $w_h\in V_h$. The iterate is then updated as
$u_{k+1,h}=u_{k,h}+\alpha_k w_h$ with step size $0<\alpha_k\leq 1$ such that under standard assumptions the iteration converges 
quadratically. The full implementation, including the assembly of the
resulting linear systems, is documented in the accompanying Zenodo
repository~\cite{sourcepackage}.

\section{Nonlinear constraints}
\label{app:nonlinear}

We explore the suitability of the method for nonlinear constraints by computing the two-membrane contact problem from Section \ref{sec:2membrane} with the nonlinear constraint $u_1 - u_2 \leq g(1 + 0.2 u_1^2)(1+0.2 u_2^2)$ in $\Omega$ instead of the constraint $u_1 - u_2 \leq g$. This results in the following Nitsche terms:
\begin{align*}
	\beta(u_h) & = u_{h,2} - u_{h,1} + g (1 + 0.2 u_{h,1}^2)(1 + 0.2 u_{h,2}^2), \\
	\lambda(u_h) & = (\kappa_1\Delta_h u_{h,1} + f_1)\cdot 0.4\cdot g \cdot u_{h,1}, \\
	\gamma(h,\kappa) & = \frac{\alpha h^2}{\kappa_1}.
\end{align*}

\begin{figure}[H]
	\begin{center}
		\includegraphics[width=0.95\textwidth]{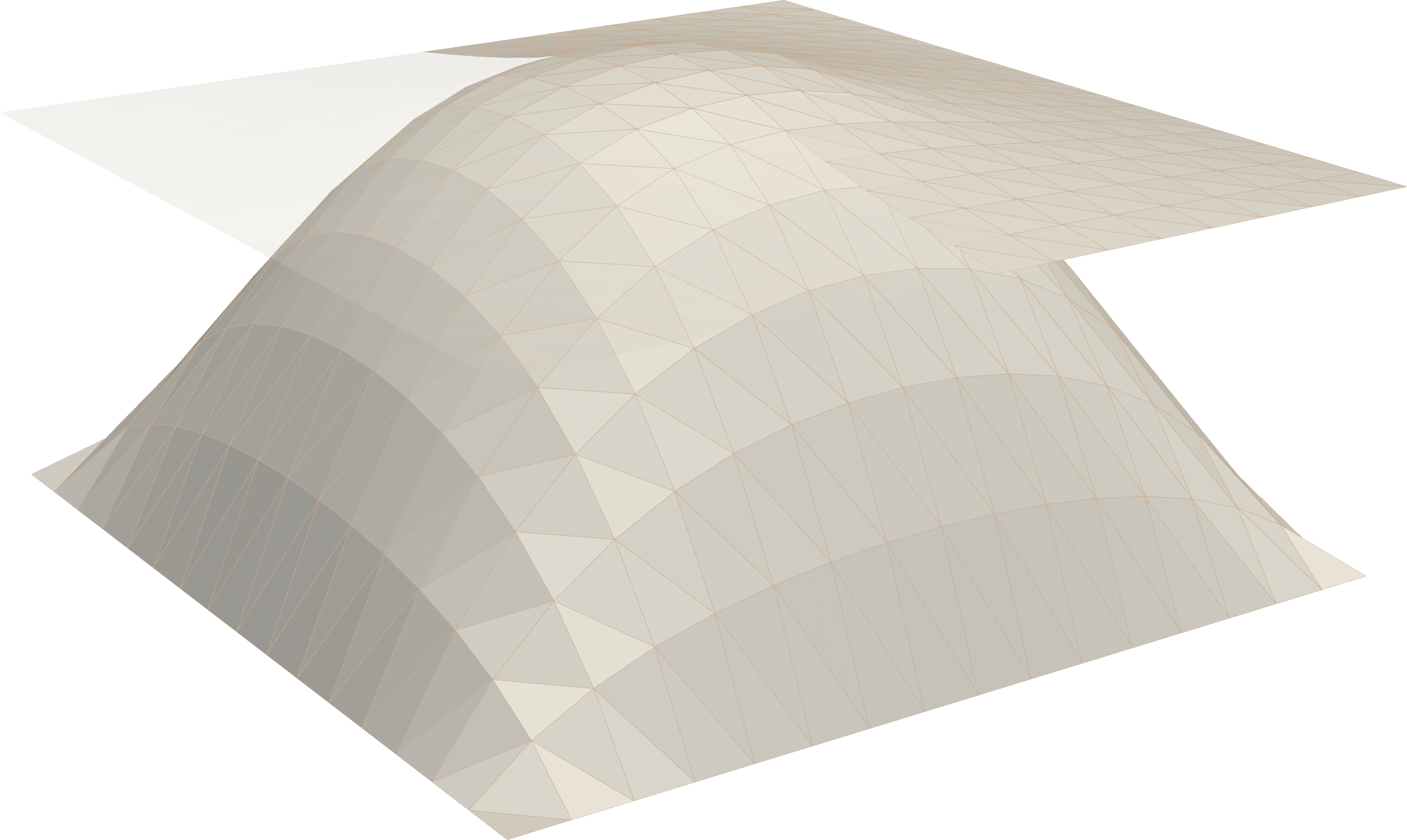}
	\end{center}
	\caption{Numerical solution for two membranes in contact from Section \ref{sec:2membrane} with the nonlinear constraint $u_1 - u_2 \leq g(1 + 0.2 u_1^2)(1+0.2 u_2^2)$ in $\Omega$.}\label{fig:2obstaclenonlinear}
\end{figure}

The problem was solved numerically over $\Omega=(0,1)^2$ with the parameterization $g=0.05$, $f_1=1, f_2=0, \kappa_1 = \kappa_2 = 1$ and 
$\alpha=10^{-2}$ 
using linear triangular elements. The numerical solution is displayed in Figure \ref{fig:2obstaclenonlinear}.

\bibliographystyle{elsarticle-num}
\bibliography{refs}

\end{document}